\documentclass{article} 
\usepackage{amssymb,latexsym}

\newtheorem{thm}{Theorem}[section]

\newtheorem{lem}[thm]{Lemma}

\newtheorem{prop}[thm]{Proposition}
\newtheorem{defn}[thm]{Definition}

\title{\bf\Large Moduli space of symmetric connections}

\author{ Stanislav Dubrovskiy } 

\date{October 31, 2001}

\begin{document}

\maketitle

\begin{abstract}
The 
action of origin-preserving diffeomorphisms on a space of 
jets of symmetric connections is considered. Dimensions of moduli spaces of generic connections
are calculated. Poincar\'e series of the geometric structure of connection is 
constructed, and shown to be a rational function. 
\end{abstract}



\section{Introduction}
A problem of finding functional moduli or at least establishing their finitness in various 
local differential-geometric settings was discussed by Arnold in \cite{Arn}.
Here we are interested in local differential invariants of
a geometric structure consisting of a symmetric connection, under smooth coordinate changes.
The structure of the resulting generic moduli space is reflected in the Poincar\'e series which we 
explicitly calculate, cf.(\ref{eq:P'series}). This series turns out to be a rational function, 
cf.(\ref{P'ratio}), 
indicating a finite number of invariants. 
This confirms the finitness assertion 
of Tresse \cite{T}, which he stated for any "natural" differential-geometric structure. 
Similar results for Riemannian, K\"{a}hler and Hyper-K\"{a}hler structures were 
obtained in \cite{Shmel:struc}, and an explicit normal form for 
Riemannian structure - in \cite{Shmel:mod}, by Shmelev. Earlier Vershik and Gershkovich investigated 
jet asymptotic dimension of moduli spaces of jets of generic distributions at $0\in\mathbb{R}^{n}$ 
in \cite{Versh:Gersh}, and their normal form in \cite{Gersh}.

I would like to acknowledge guidance and support for this work from M.Shubin. I am 
grateful to C.L.Terng, A.Vershik and J.Weyman for fruitful discussions.


\section{Preliminaries and main result}
Let $\mathcal{F}$ and $\mathcal{F}_{k}$ be the spaces of germs and $k$-jets  
of symmetric $C^{\infty}$-connections at a point on $\mathbb{R}^{n}$. 
From now on all connections we consider are assumed 
symmetric. As usual, two $C^{\infty}$-functions on 
$\mathbb{R}^{n}$ have the same $k$-jet at a point if their first $k$ derivatives are equal in any local coordinates.
We say that two connections $\nabla$ and $\tilde{\nabla}$ have the same $k$-jet at 0 if for any two 
$C^{\infty}$-vector fields $X$,$Y$ and any $C^{\infty}$-function $f$, the functions $\nabla_{X}Y(f)$
and $\tilde{\nabla}_{X}Y(f)$ have the same $k$-jet at $0$. ( This is equivalent to Cristoffel symbols of
$\nabla$ and $\tilde{\nabla}$ having the same $k$-jet )\\
We will frequently denote connection and its Cristoffel symbol with the same letter,
e.g. $\Gamma$ ; $j^{k}\Gamma$ would stand for its $k$-jet .

There is an action of the group of germs of origin-preserving diffeomorphisms 
$G:=\mathrm{Diff}(\mathbb{R}^{n},0)$ on $\mathcal{F}$ and $\mathcal{F}_{k}$.
For $\varphi \in G$, $\nabla$(or $\Gamma$) $\in \mathcal{F}$ and $j^{k}\Gamma \in \mathcal{F}_{k}$:
$$\Gamma \mapsto \varphi^{*}\Gamma\,, \quad j^{k}\Gamma \mapsto j^{k}(\varphi^{*}\Gamma)\,,$$
where $${( \varphi^{*}\nabla )}_{X} Y={\varphi}_{*}^{-1}( \nabla_{  \phantom{|}\atop 
{\varphi}_{*}X  } {\varphi}_{*}Y )$$
Let us introduce a filtration of $G$ by normal subgroups:
$$ G=G_{1}\rhd G_{2}\rhd G_{3}\rhd \ldots,$$
where $$G_{k}=\{
\,\varphi \in G\ |\ \varphi(x)=x+(\varphi_{1}(x),\ldots\varphi_{n}(x)),\ \varphi_{i}=O(|x|^{k}),
\,i=1,\ldots,n\,
\}
$$
The subgroup $G_{k}$ acts trivially on $\mathcal{F}_{p}$ for $k\geq p+3$.
It means that the action of $G$ coincides with that of $G/G_{p+3}$ on each $\mathcal{F}_{p}$.
Now $G/G_{p+3}$ is a finite-dimensional Lie group, which we will call $K_{p}$.
$j^{k}\Gamma\in\mathcal{F}_k$
Denote by $\mathrm{Vect}_{0}(\mathbb{R}^{n})$ the Lie algebra of $C^{\infty}$-vector fields, 
vanishing at the origin.
It acts on $\mathcal{F}
$ as follows:
\begin{defn}
For $V\in\mathrm{Vect}_{0}(\mathbb{R}^{n})$ 
generating a local 1-parameter subgroup $g^{t}$ of $
\mathrm{Diff}(\mathbb{R}^{n},0)$,
\emph{the Lie derivative} of a connection $\nabla$ in the direction $V$ is a (1,2)-tensor:
$$
\mathcal{L}_{V}\nabla
=\left. \frac{d}{dt}\right|_{t=0} 
{g^{t}}^{*}\nabla
$$
\end{defn}
\begin{lem} 
\begin{equation}
\label{eq:LD}
(\mathcal{L}_{V}\nabla)(X,Y)=[V,\nabla_{X}Y]-\nabla_{[V,X]}Y-\nabla_{X}[V,Y]  
\end{equation}
\end{lem} 
{\bf Proof}\quad
Below the composition $\circ$ is understood as that of differential operators acting on functions.
$$
(\mathcal{L}_{V}\nabla)(X,Y)=\left. \frac{d}{dt}\right|_{t=0}  
g_{*}^{-t}[ \nabla_{  \phantom{|}\atop { g^{t}_{*}X } } g^{t}_{*}Y  ]=
\left. \frac{d}{dt}\right|_{t=0}
\left[\,(g^{t})^{*}\circ [\nabla_{  \phantom{|}\atop { g^{t}_{*}X } } g^{t}_{*}Y ] \circ (g^{-t})^{*}\right]=
$$
$$
\left. \frac{d}{dt}\right|_{t=0} (g^{t})^{*} \circ \nabla_{X}Y +
\nabla_{X}Y \circ \left. \frac{d}{dt}\right|_{t=0} (g^{-t})^{*} +
\nabla_{ \left. \frac{d}{dt}\right|_{t=0}g^{t}_{*}X }Y+
\nabla_{X}\left. \frac{d}{dt}\right|_{t=0}g^{t}_{*}Y =
$$
$$
V\circ \nabla_{X}Y - \nabla_{X}Y \circ V - \nabla_{\left. \frac{d}{dt}\right|_{t=0}g_{*}^{-t}X}Y
 - \nabla_{X}\left. \frac{d}{dt}\right|_{t=0}g_{*}^{-t}Y =
$$
$$
\mathcal{L}_{V}(\nabla_{X}Y)-\nabla_{\mathcal{L}_{V}X}Y - \nabla_{X}(\mathcal{L}_{V}Y)
$$
    \hfill$\Box$\\
This defines the action on the germs of connections.
Now we can define the action of $\mathrm{Vect}_{0}(\mathbb{R}^{n})$ on jets $\mathcal{F}_{k}$.
For $V\in \mathrm{Vect}_{0}(\mathbb{R}^{n})$:
$$\mathcal{L}_{V}(j^{k}\Gamma)=j^{k}(\mathcal{L}_{V}\Gamma)\ ,$$ 
where $\Gamma$ on the right 
is an arbitrary representative of the $j^{k}\Gamma$ on the left.\\
This is well-defined, since in the coordinate version of (\ref{eq:LD}):
\begin{equation}
\label{eq:LDcoor}
(\mathcal{L}_{V}\Gamma)_{ij}^{l}=
V^{k}\frac{\partial \Gamma_{ij}^{l}}{\partial x^{k}}-
\Gamma_{ij}^{k}\frac{\partial V^{l}}{\partial x^{k}}+
\Gamma_{kj}^{l}\frac{\partial V^{k}}{\partial x^{i}}+
\Gamma_{ik}^{l}\frac{\partial V^{k}}{\partial x^{j}}+
\frac{\partial^{2}V^{l}}{\partial x^{i} \partial x^{j}}
\end{equation}
elements of $k$-th order and less are only coming from $j^{k}\Gamma$, because $V(0)=0$.
Einstein summation convention in (\ref{eq:LDcoor}) above and further on is assumed.\\ 
Consequently, the action is invariantly defined.
This can also be expressed as commutativity of the following diagram:
$$
\begin{array}{ccccccccccccc}
j^{0}\mathcal{F}& \longleftarrow & \ldots & \longleftarrow & j^{k-1}\mathcal{F} & \stackrel{\pi_{k}}{\longleftarrow} & j^{k}\mathcal{F} & 
\longleftarrow & \ldots & \longleftarrow & \mathcal{F} &  &  
\\
\downarrow\lefteqn{\mathcal{L}_{V}}& & & & \downarrow\lefteqn{\mathcal{L}_{V}} & & \downarrow\lefteqn{\mathcal{L}_{V}} & & 
 & & \downarrow\lefteqn{\mathcal{L}_{V}}& &\!\!\!, 
\\
j^{0}\Pi & \longleftarrow & \ldots & \longleftarrow & j^{k-1}\Pi & \stackrel{\pi_{k}}{\longleftarrow} & j^{k}\Pi  & 
\longleftarrow & \ldots & \longleftarrow & \Pi & & 
\\
\end{array}
$$
where $\pi_{k}$ is projection from $k$-jets onto $(k-1)$-jets, 
$\mathcal{F}$ and $\Pi$
denote spaces of germs of connections and that of (1,2)-tensors
respectively, at $0$.

\emph{Poincar\'e series} will encode information about these actions for all $k$.
The space $\mathcal{M}=\mathcal{F}/\mathrm{Diff}(\mathbb{R}^{n},0)$ of 
$\mathrm{Diff}(\mathbb{R}^{n},0)$-orbits on $\mathcal{F}$ is called the \emph{moduli space} 
of connections at $0$ on $\mathbb{R}^{n}$. We do not introduce any topology on $\mathcal{M}$. 
Similarly, the orbit space 
$\mathcal{M}_{k}=\mathcal{F}_k/\mathrm{Diff}(\mathbb{R}^{n},0)=\mathcal{F}_k/K_{k}$ 
is called the moduli space of connection $k$-jets. 

The action of $K_k$ is algebraic, a subspace $\mathcal{F}_{k}^{0}\subset\mathcal{F}_{k}$ of points on generic orbits
(those of largest dimension) is a smooth manifold, open and dense in $\mathcal{F}_{k}$. 
Subspace of points on orbits of any other given dimension is a manifold as well, 
albeit of a lesser dimension. 
We could consider the $G$-quotient for each of those subspaces, and have a moduli space of its own for
each of the orbit types. Let $\mathcal{O}_{k}$ denote a generic orbit. Denote by $\mathcal{M}_{k}^{0}$ 
the moduli space of generic connections: 
$$\mathcal{M}_{k}^{0}=\mathcal{F}^0_k/\mathrm{Diff}(\mathbb{R}^{n},0)=\mathcal{F}^0_k/K_{k}$$
or generic subspace of moduli space ${F}^0_k$ .
Its dimension is found as:
\begin{equation}
\label{eq:dim}
\dim\mathcal{M}_{k}^{0}=\dim\mathcal{F}_{k}^{0}-\dim\mathcal{O}_{k}
\end{equation} 
Unfortunately it is not true that the generic subspace has the maximal dimension for general Lie group actions, 
even those that are algebraic, for an explicit counterexample 
see \cite{Weym}.
Thus we define $$\dim\mathcal{M}_{k}=\dim\mathcal{M}_{k}^{0}
$$ for mere simplicity of notation.
One more piece of notation:  
\begin{displaymath}
a_{k}=\left\{ \begin{array}{ll} \dim \mathcal{M}_{k} & ,k=0 \\ 
                                \dim \mathcal{M}_{k}- \dim \mathcal{M}_{k-1} & ,k\geq 1
\end{array} \right. 
\end{displaymath}
and we can introduce our main object of interest.

\begin{defn}
The formal power series
$$p_{\Gamma}(t)=\sum_{k=0}^{\infty}a_{k}t^{k}$$
is called the \emph{Poincar\'e series} for the moduli space $\mathcal{M}$. 
\end{defn}
Our main result is the following
\begin{thm}
\label{thm:series}
Poncar\'e series coefficients $a_{k}=a(k)$ are polynomial in $k$, and the series has the form:\\\\
$p_{\Gamma}(t)=(t-t^2)\delta_{2}^{n}+
{\displaystyle n\sum_{k=1}^{\infty}\;\biggl[\frac{n(n+1)}{2}{n+k-1 \choose n-1}
-n{n+k+1 \choose n-1}\biggr]t^{k}}
$\\\\
$\hspace*{22em}(\ \delta\textrm{\emph{ is Kronecker symbol }}).$\\\\
It represents a rational function.
\end{thm}
{\bf Remark}\ \ This complies with Tresse' 
assertion that algebras of "natural" differential-geometric structures be finitely-generated.\\\\
{\bf Proof}\ \ of this theorem is relegated to section \ref{sec:Proof}.\\\\
To explain significance of Poincar\'e series represented by a rational function, we make the following:\\
{\bf Remark}\ \ 
If a geometric structure is described by a finite number of functional moduli, then its Poncar\'e series is rational.
In particular, if there are $m$ functional invariants in $n$ variables, then
$$p(t)=\frac{m}{(1-t)^n}$$
Indeed, dimension of moduli spaces of $k$-jets is just the number of monomials up to the order $k$
in the formal power series of the $m$ given invariants:$$\dim \mathcal{M}_{k}=m{n+k \choose n}$$
For more details and slightly more general formulation see {\bf Theorem 2.1} in \cite{Shmel:mod}.
\section{Stabilizer of a generic k-jet}
In order to calculate Poincar\'e series, we need to find $\dim\mathcal{O}_{k}=\dim K_{k}-\dim G_{\Gamma}$, 
and hence, the size of a stabilizer $G_{\Gamma}$ for a generic connection $\Gamma$.

Let us find the subalgebra generating $G_{\Gamma}$ - the stabilizer of a $k$-jet $\Gamma$.\\
It consists of such $V\in\mathrm{Vect}_{0}(\mathbb{R}^{n})$, that:
\begin{equation}\label{eq:stab}
\mathcal{L}_{V}(j^{k}\Gamma)=0
\end{equation}
We will argue in local coordinates. Let us introduce grading in homogeneous components:
$$V=V_{1}+V_{2}+\ldots$$ 
\begin{center} ( $V_{0}=0$, since we require $V$ to preserve the origin ) , 
\end{center}
$$\Gamma=\Gamma_{0}+\Gamma_{1}+\ldots$$-- analogous grading on the connection $\Gamma$.\\
Note that each $\Gamma_{k}$ has same symmetry on indexes as the original $\Gamma$.\\ 
Using these decompositions we can rewrite (\ref{eq:stab}) as follows:
$$
\mathcal{L}_{V}(j^{k}\Gamma)=j^{k}\mathcal{L}_{V}(\Gamma)= 
j^{k}\mathcal{L}_{V_{1}+V_{2}+\ldots}(\Gamma_{0}+\Gamma_{1}+\ldots+\Gamma_{k}+\ldots)=
$$
$
=\underbrace{\mathcal{L}_{V_{1}}\Gamma_{0}+ 
\frac{\partial^{2}V_{2}}{\partial x^{2}}}_{\textrm{\footnotesize{0th order}}}+  
\underbrace{ \mathcal{L}_{V_{1}}\Gamma_{1}+{\tilde\mathcal{L}}_{V_{2}}\Gamma_{0}+ 
\frac{\partial^{2}V_{3}}{\partial x^{2}} }_{\textrm{\footnotesize{1st order}}}+ 
\ldots
$ 

$\qquad\qquad\qquad\qquad\qquad
\ldots+ 
\underbrace{\tilde{\mathcal{L}}_{V_{k+1}}\Gamma_{0}+
\tilde{\mathcal{L}}_{V_{k}}\Gamma_{1}+\ldots+\mathcal{L}_{V_{1}}\Gamma_{k}+
\frac{\partial^{2}V_{k+2}}{\partial x^{2}} }_
{\textrm{\footnotesize{k-th order}}}  
\quad,
$
$$\mathrm{where}
\quad\left(\frac{\partial^{2}V_{2}}{\partial x^{2}}\right)^{l}_{ij}=
\frac{\partial^{2}V_{2}^{l}}{\partial x^{i}x^{j}}$$
and
$$\quad\tilde{\mathcal{L}}_{V}\Gamma=
\mathcal{L}_{V}\Gamma-\frac{\partial^{2}V}{\partial x^{2}}\,,$$ the same with indexes:
$$(\tilde{\mathcal{L}}_{V}\Gamma)_{ij}^l=
V^{k}\frac{\partial \Gamma_{ij}^{l}}{\partial x^{k}}-
\Gamma_{ij}^{k}\frac{\partial V^{l}}{\partial x^{k}}+
\Gamma_{kj}^{l}\frac{\partial V^{k}}{\partial x^{i}}+
\Gamma_{ik}^{l}\frac{\partial V^{k}}{\partial x^{j}}$$
are just the first 4 terms from (\ref{eq:LDcoor}) of $({\mathcal{L}}_{V}\Gamma)_{ij}^l$.  
The stabilizer condition therefore results in a system:
\begin{equation}\label{sys:SYS}
\left\{
\begin{array}{l} 
\mathcal{L}_{V_{1}}\Gamma_{0}+{\displaystyle\frac{\partial^{2}V_{2}}{\partial x^{2}}} =  0 \\
\mathcal{L}_{V_{1}}\Gamma_{1}+{\tilde\mathcal{L}}_{V_{2}}\Gamma_{0}+
{\displaystyle\frac{\partial^{2}V_{3}}{\partial x^{2}}}  =  0\\
\qquad \vdots  \\
\tilde{\mathcal{L}}_{V_{k+1}}\Gamma_{0}+
\tilde{\mathcal{L}}_{V_{k}}\Gamma_{1}+\ldots+\mathcal{L}_{V_{1}}\Gamma_{k}+
{\displaystyle\frac{\partial^{2}V_{k+2}}{\partial x^{2}}}  =  0
\end{array}
\right. 
\end{equation}
Our present goal is finding all $(V_{1},V_{2},\ldots,V_{k+2})$ solving the above system for a generic $\Gamma$.
Assume $V_{1}$ is arbitrary, to find $V_{2}$ from the first equation we need the following
\begin{prop}
\label{lem:mixder} 
Given a family $\{f_{ij}\}_{1\leq i,j\leq n}$ of smooth functions,
solution $u$ for the system:
$$
\left\{
\begin{array}{l}
u_{,kl}=f_{kl}\\
1\leq k,l\leq n
\end{array}\right.
$$
(indexes after a comma henceforth will denote differentiations in corresponding variables)
exists if and only if
\begin{equation}
\label{eq:compatibility}
\left\{
\begin{array}{l}
f_{ij}=f_{ji}\\
f_{ij,k}=f_{kj,i}
\end{array}\right.
\end{equation}
If $f_{ij}$ are homogeneous polynomials of degree $s\ge0$, then $u$ can be uniquely chosen as a polynomial of degree $s+2$.
\end{prop}
{\bf Proof}\ \ is a straightforward integration of the right-hand sides.
\hfill$\Box$\\
Therefore, if we treat highest-order $V_{k}$ in each equation of 
(\ref{sys:SYS}) as an unknown, we see that various (combinations 
of) $\mathcal{L}_{V}\Gamma$ must satisfy (\ref{eq:compatibility}).
The first condition is satisfied automatically. The second one gives:
$$(\mathcal{L}_{V}\Gamma)_{ij,p}^{l}=(\mathcal{L}_{V}\Gamma)_{pj,i}^{l}$$
This condition  for the 1st equation in (\ref{sys:SYS}) is satisfied trivially, 
since $V_{1}$ is of the 1st degree, and $\Gamma_{0}$ - $0^{th}$ degree in $x$.
Hence, $V_{2}$ exists, and since it must be of the second degree, is unique.
Let us now find $V_3$ from the next equation:
$$\mathcal{L}_{V_1}\Gamma_1 + \widetilde{\mathcal{L}}_{V_2}\Gamma_0 + 
\frac{\partial V_3}{\partial x^{2}} = 0$$
It follows from the Proposition\,\ref{lem:mixder}\,, that for existence of $V_3$ it is necessary 
( and sufficient ) to have the following condition:
\begin{equation}
\label{eq:9}
(\mathcal{L}_{V_1}\Gamma_1 + \widetilde{\mathcal{L}}_{V_2}\Gamma_0)_{ij,p}=
(\mathcal{L}_{V_1}\Gamma_1 + \widetilde{\mathcal{L}}_{V_2}\Gamma_0)_{pj,i}
\end{equation}
We assert that, except for cases $n\le 2$ considered in the next section, it is not satisfied for a generic connection unless $V_1=0$ .
In other words (\ref{eq:9}), considered as a condition on $V_1$ implies $V_1=0$ 
( and hence $V_2=V_3=...=0$ too ). Indeed, write: 
\begin{equation}
\label{V:comp}
V_{1}^{k}=\sum_{s=1}^{n}b_{s}^{k}x^{s}
\end{equation}
and let us consider (\ref{eq:9}) as a homogeneous linear system on components $\{b^{k}_{s}\}$ of $V_1$. 
In what follows we will see that for a generic connection this system is nondegenerate.
More precisely, we will find a connection, and a suitable minor of the system
( since there are more equations than variables, we will choose a subset of equations
to obtain a square matrix ), and show that it is non-degenerate.
Since it is an open condition, it would be generically true. 
We will construct connection and minor step by step, trying to obtain one
as close to diagonal as possible.
Let us look for such connection $\Gamma$ among those with $\Gamma_0$=0.
Then we see (\ref{eq:9}) shrink to: 
$$(\mathcal{L}_{V_1}\Gamma_1)_{ij,p}^{l}=(\mathcal{L}_{V_1}\Gamma_1 )_{pj,i}^{l}$$
Let us expand it using (\ref{eq:LDcoor}):
$$V_{1 ,p}^{k}\Gamma_{1 ij,k}^{l}-V_{1 ,k}^{l}\Gamma_{1 ij,p}^{k}+
V_{1 ,i}^{k}\Gamma_{1 kj,p}^{l}+V_{1 ,j}^{k}\Gamma_{1 ik,p}^{l}= 
(i \leftrightarrow p)\quad,$$
where the right-hand side is the same as the left-hand side with the two indexes swapped.
After simplifying we get:
\\[3.5ex]
$(\Gamma_{1 ij,k}^{l}-\Gamma_{1 kj,i}^{l})V_{1 ,p}^{k}+
(\Gamma_{1 kj,p}^{l}-\Gamma_{1 pj,k}^{l})V_{1 ,i}^{k}+$
\begin{equation}
\label{eq:unGEN}
\qquad
(\Gamma_{1 pj,i}^{k}-\Gamma_{1 ij,p}^{k})V_{1 ,k}^{l}+
(\Gamma_{1 ik,p}^{l}-\Gamma_{1 pk,i}^{l})V_{1 ,j}^{k}=0
\end{equation}
\vspace{2ex}\\
where $1\le i<p\leq n$, and $j,l$ vary from 1 to $n$.\\ 
That makes for $n^{2}\!\cdot\!\!{\displaystyle \frac{n(n-1)}{2}}$ equations 
on $n^{2}$ variables.\\
Also write:
\[\Gamma_{1 ij}^{l}=\sum_{m=1}^{n}c_{ij}^{lm}x^{m}, 
\qquad c_{ij}^{lm}=c_{ji}^{lm}\textrm{  (recall that the connection is assumed symmetric) }\]
Then, recalling also (\ref{V:comp}), (\ref{eq:unGEN}) becomes:
\begin{equation}\label{eq:sys}
(c_{ij}^{lk}-c_{kj}^{li})b_{p}^{k}+
(c_{kj}^{lp}-c_{pj}^{lk})b_{i}^{k}+
(c_{pj}^{ki}-c_{ij}^{kp})b_{k}^{l}+
(c_{ik}^{lp}-c_{pk}^{li})b_{j}^{k}=0
\end{equation}
Recall that summation over repeated indexes above is assumed.\\
Let us now impose a further restriction on $\Gamma$, namely that 
\begin{equation}
\label{ex:condition}
c_{ij}^{lp}\neq0\textrm{  only if  }\{i,j,l,p\}=\{\alpha,\beta\},\alpha\neq\beta\ 
\end{equation}
In other words nonzero coefficients may only occur among those with indexing set consisting of two distinct numbers, and must be zero otherwise.
\begin{lem}The homogeneous linear system (\ref{eq:sys}) with coefficients 
as restricted in (\ref{ex:condition}) has a $n^{2}\times n^{2}$
nondegenerate minor. 
\end{lem}
{\bf Proof} 
Let us specify the minor by letting $j$ and $l$ in the index set $\{i,j,l,p\}$  be arbitrary. 
Accordingly, we will be labeling equations in the system by this pair of indexes $(jl)$, in lexicographic order.
The remaining two indexes $i$ and $p$ will be determined by \{$j,l$\} in the following manner:\\
i) if $j<l$ set $i=j$, $p=l$\\
Equation ($jl$) becomes then:
$$(jl)\qquad (c_{jj}^{lk}-c_{kj}^{lj})b_{l}^{k}+ 
(c_{kj}^{ll}-c_{lj}^{lk})b_{j}^{k}+ 
(c_{lj}^{kj}-c_{jj}^{kl})b_{k}^{l}+ 
(c_{jk}^{ll}-c_{lk}^{lj})b_{j}^{k}=0$$
Diagonal coefficients ( at $b_{j}^{l}$ ) are:
$$(c_{jl}^{jj}-c_{jj}^{jl})+(c_{lj}^{ll}-c_{ll}^{lj})=:A_{jl}=A_{lj}$$
Using (\ref{ex:condition}), the only non-diagonal terms that remain are:
$$(c_{jj}^{ll}-c_{lj}^{lj})b_{l}^{l}+ 
(c_{jj}^{lj}-c_{jj}^{lj})b_{l}^{j}+ 
(c_{jj}^{ll}-c_{lj}^{lj})b_{j}^{j}+ 
(c_{lj}^{lj}-c_{jj}^{ll})b_{l}^{l}+
(c_{jj}^{ll}-c_{lj}^{lj})b_{j}^{j}=
2(c_{jj}^{ll}-c_{lj}^{lj})b_{j}^{j}$$
To eliminate them we need to set $c_{jj}^{ll}=c_{lj}^{lj}$ for $j<l$ .\\\\
ii) if $j=l<n$ set $i=l$, $p=l+1$ \\
Then the corresponding equations have the form:
$$\hspace{-1em}(ll)\ (c_{ll}^{lk}-c_{kl}^{ll})b_{l+1}^{k}+ 
(c_{kl}^{l(l+1)}-c_{(l+1)l}^{lk})b_{l}^{k}+ 
(c_{(l+1)l}^{kl}-c_{ll}^{k(l+1)})b_{k}^{l}+ 
(c_{lk}^{l(l+1)}-c_{(l+1)k}^{ll})b_{l}^{k}=0$$
with diagonal coefficients ( at $b_{l}^{l}$ ):
$(c_{ll}^{l(l+1)}-c_{(l+1)l}^{ll})\;,$\\
and suspicious off-diagonal: \\
$$(c_{ll}^{ll}-c_{ll}^{ll})b_{l+1}^{l}+ 
(c_{ll}^{l(l+1)}-c_{(l+1)l}^{ll})b_{l+1}^{l+1}+ 
(c_{(l+1)l}^{l(l+1)}-c_{(l+1)l}^{l(l+1)})b_{l}^{l+1}+$$
$$+(c_{(l+1)l}^{(l+1)l}-c_{ll}^{(l+1)(l+1)})b_{l+1}^{l}+ 
(c_{l(l+1)}^{l(l+1)}-c_{(l+1)(l+1)}^{ll})b_{l}^{l+1}=$$
$$=(c_{ll}^{l(l+1)}-c_{(l+1)l}^{ll})b_{l+1}^{l+1}+ 
(c_{(l+1)l}^{l(l+1)}-c_{(l+1)(l+1)}^{ll})b_{l}^{l+1}+
(c_{(l+1)l}^{(l+1)l}-c_{ll}^{(l+1)(l+1)})b_{l+1}^{l}$$
Hence to nullify them, we could require that:
$$c_{ll}^{l(l+1)}=c_{(l+1)l}^{ll},
c_{(l+1)l}^{(l+1)l}=c_{ll}^{(l+1)(l+1)}, 
c_{l(l+1)}^{l(l+1)}=c_{(l+1)(l+1)}^{ll}\textrm{ for any }l: 1\leq l\leq n-1$$ 
iii) if $j>l$ set $i=l$, $p=j$\\
These are the equations:
$$(jl) \qquad (c_{lj}^{lk}-c_{kj}^{ll})b_{j}^{k}+ 
(c_{kj}^{lj}-c_{jj}^{lk})b_{l}^{k}+ 
(c_{jj}^{kl}-c_{lj}^{kj})b_{k}^{l}+ 
(c_{lk}^{lj}-c_{jk}^{ll})b_{j}^{k}=0$$
The diagonal coefficients ( at $b_{j}^{l}$ ) are: 
$$(c_{jj}^{jl}-c_{jl}^{jj})+(c_{ll}^{lj}-c_{lj}^{ll})=-A_{jl}$$
Off-diagonal coefficients are: $2(c_{lj}^{lj}-c_{jj}^{ll})b_{j}^{j}$,\\
hence it's necessary to require that 
$$c_{lj}^{lj}=c_{jj}^{ll}\textrm{ for }j>l$$
iv) finally, if $j=l=n$ set:
$$\left\{
\begin{array}{l}
i=1, p=n \textrm{ for $n$-odd}
\\
i=2, p=n \textrm{ for $n$-even.}
\end{array}\right.$$
We have for $n$-odd :
$$(nn)
\qquad (c_{1n}^{nk}-c_{kn}^{n1})b_{n}^{k}+ 
(c_{kn}^{nn}-c_{nn}^{nk})b_{1}^{k}+ 
(c_{nn}^{k1}-c_{1n}^{kn})b_{k}^{n}+ 
(c_{1k}^{nn}-c_{nk}^{n1})b_{n}^{k}=0$$
The diagonal coefficients ( at $b_{n}^{n}$ ) are: $c_{1n}^{nn}-c_{nn}^{n1}$,\\
the off-diagonal ones are: 
$$(c_{1n}^{nn}-c_{nn}^{n1})\textrm{ at }b_{1}^{1}, 
(c_{nn}^{11}-c_{1n}^{1n})\textrm{ at }b_{1}^{n}, 
(c_{11}^{nn}-c_{n1}^{n1})\textrm{ at }b_{n}^{1}$$
Hence, we would need: $$c_{1n}^{nn}=c_{nn}^{n1}, 
c_{nn}^{11}=c_{1n}^{1n}, 
c_{11}^{nn}=c_{n1}^{n1}$$
to get rid of them.\\
For $n$-even just exchange index 1 for index 2 in all equations throughout iv) above.\\

From this it's clear that all equations can be "diagonalized", except those labeled $(ll)$
( the requirements above are contradictory in ii) and iv) ). Each of these equations then 
will have a single off-diagonal coefficient. \\
Summarizing, we set:
$$
\begin{array}{l}
c_{jj}^{ll}=c_{lj}^{lj}, \\
c_{jj}^{jl}-c_{jl}^{jj}=1,\quad j\neq l. 
\end{array}
$$
That implies $A_{jl}=-2$ for $j\neq l$,\\
for $l<n$ in $(ll)^{\mathrm{th}}$ equation there is a unit off-diagonal coefficient at $b_{l+1}^{l+1}$, and\\
for $l=n$ in $(nn)^{\mathrm{th}}$ equation there is a negative unit off-diagonal coefficient at $b_{1}^{1}$ 
for $n$-odd, or at $b_{2}^{2}$ for $n$-even.\\
Now we switch the order of equations and variables , so that $(kk)^{\mathrm{th}}$ 
equations and variables appear first, followed by all the rest in the preset lexicographic order.\\
Then the minor for $n$-odd will assume the form:
$$\left(
\begin{array}{cccccccc}
        1       &	        1            &                    &                 &        \\
                &\hspace{-2.5em}\ddots &\hspace{-4em}\ddots &                 &       \\
                &\qquad \ddots         &\hspace{1.4em}\ddots&                 & 0& & &\\
                &  			   &			      &                 &   \\
                &                      &\hspace{1em}1       &\hspace{-2.15em}1&   \\
                &		               &				&			&         \\
\hspace{-.7em}-1& 		         &  	            &\hspace{-2.7em}-1& &      &    \\
		    &			         &		            &			&	\\
                &                      &                    &\hspace{3em} -2  &          &      &    \\
                &        	         & \hspace{-3.3em}0   &                 &\hspace{.2em}\ddots&      & 	\\
                &       		   &                    &                 &	              &\hspace{-.4em}\pm2\\
                &		               &                    &                 &	     &        &\hspace{-.2em}\ddots\\
                &                      &                    &                 &          &        &        &2  \\ 
\end{array} \right)$$\\
For $n$-even, just relocate $-1$ in the $n$-th row from the first position to the second.
The minor above is easily seen to be nondegenerate, as required. \hfill$\Box$\\\\
This argument proves 
\begin{prop}\label{stab:highdim}
The stabilizer of a k-jet of a generic connection for $n\geq 3$ is:
$G_{1}/G_{2}$ 
 for $k=0$
, and $0$ 
 for $k\geq 1$ .
\end{prop}
\section{Exceptions: the stabilizer in low dimensions}
Let us start with the case $n=1$ .
In this case any index can assume only a single possible value: $1$ . 
Then the compatibility conditions (\ref{eq:compatibility}) are vacuous in all $k$-jets.
Hence the stabilizer for $k$-jet is determined by $V_1$ with no restrictions and is equal to $G_{1}/G_{2}$ 
for all $k$. That is a $1$-dimensional stabilizer, resulting in a $k+1$-dimensional orbit and $0$-dimensional
moduli space. Poincare series then is identically equal to $0$.

For the case $n=2$, the compatibility conditions arising from second equation in (\ref{sys:SYS}) have non-trivial solutions, that is, unlike the higher $n$'s, the stabilizer of the 1st jet is non-trivial.
The reason is that the analogue of equation (\ref{eq:9}), considered as a ( 4 by 4 ) homogeneous linear system 
on coefficients $b^l_k$ of $V_1$ is degenerate.
Compatibility conditions for the third equation however do make up a non-degenerate linear system and have only
trivial solution. Hence for $n=2$ stabilizer is trivial, starting with second jet ( for higher $n$, trivial starting with first jet, and for $n=1$, always non-trivial ).

Let us start with considering (\ref{eq:9}) for $n=2$, showing it is degenerate, and finding its rank.
Notice that we consider (\ref{eq:9}) in the most general form, 
with arbitrary $\Gamma_0$ and $\Gamma_1$ .
Since it involves $\widetilde{\mathcal{L}}_{V_2}\Gamma_0$, we need to express $V_2$ from the 
first equation of the system: 
$$\mathcal{L}_{V_{1}}\Gamma_{0}+{\displaystyle\frac{\partial^{2}V_{2}}{\partial x^{2}}} =  0 \,.$$ 
Setting $\Gamma^k_{0ij}=:\gamma^k_{ij}$, it can be rewritten in index form as:
$$V_{2,ij}^l=\gamma_{ij}^{k}b^{l}_k-\gamma_{kj}^{l}b_i^{k}-\gamma_{ik}^{l}b_j^{k}=:v^{l}_{ij}$$
Actually, second derivatives of $V_2$ is all we need in (\ref{eq:9}), where they appear in
$$(\widetilde{\mathcal{L}}_{V_2}\Gamma_0)^l_{ij,p}
=-\gamma^k_{ij}V_{2,kp}^l+\gamma^l_{kj}V_{2,ip}^k+\gamma^l_{ik}V_{2,jp}^k\ ,$$
which we can now rewrite as:\\\\
$(\widetilde{\mathcal{L}}_{V_2}\Gamma_0)^l_{ij,p}=
-\gamma^k_{ij}(\gamma_{kp}^{s}b^{l}_s-\gamma_{sp}^{l}b_k^{s}-\gamma_{ks}^{l}b_p^{s})
+\gamma^l_{kj}(\gamma_{ip}^{s}b^{k}_s-\gamma_{sp}^{k}b_i^{s}-\gamma_{is}^{k}b_p^{s})$
\begin{equation}
\label{eq:LV2Gamma0}
\hspace*{14em}
+\gamma^l_{ik}(\gamma_{jp}^{s}b^{k}_s-\gamma_{sp}^{k}b_j^{s}-\gamma_{js}^{k}b_p^{s})
\end{equation}
If we consider (\ref{eq:9}) as $S V_1=0$ - linear operator acting on $V_1$, and split the operator into two parts: $S=S(\Gamma_0)+S(\Gamma_1)$ , then (\ref{eq:LV2Gamma0}) allows us to rewrite 
$S(\Gamma_0) V_1=(\widetilde{\mathcal{L}}_{V_2}\Gamma_0)^l_{ij,p}-(\widetilde{\mathcal{L}}_{V_2}\Gamma_0)^l_{pj,i}$ 
as:\\\\
$
( \gamma^k_{pj}\gamma_{ki}^{s}- \gamma^k_{ij}\gamma_{kp}^{s} ) b^{l}_s +
( \gamma^l_{pk}\gamma_{js}^{k}- \gamma^k_{pj}\gamma_{ks}^{l} ) b^{s}_i +$
\begin{equation}
\label{eq:comp:gamma0}
\hspace{9em}+( \gamma^l_{pk}\gamma_{si}^{k}- \gamma^l_{ik}\gamma_{sp}^{k} ) b^{s}_j+
( \gamma^k_{ij}\gamma_{ks}^{l}- \gamma^l_{ik}\gamma_{js}^{k} ) b^{s}_p 
\end{equation}
Recall that $n=2$, and indexes $1=i<p=2$ must therefore stay fixed at $i=1$, $p=2$, 
while the remaining pair of indexes take any values.
That turns (\ref{eq:comp:gamma0}) into a system of 4 expressions, which we index, 
as in the previous section with $(j,l)$ on left:\\\\
$(11)\hspace{2em}( \gamma^1_{22}\gamma_{11}^{2}- \gamma^1_{12}\gamma_{12}^{2} ) b^{1}_1 +
		     ( \gamma^2_{11}\gamma_{22}^{1}- \gamma^1_{12}\gamma_{12}^{2} ) b^{2}_2 + $\\\\
$\hspace*{17em}    +( \gamma^1_{2k}\gamma_{12}^{k}- \gamma^1_{1k}\gamma_{22}^{k} ) b^{2}_1 +
		     ( \gamma^k_{21}\gamma_{k1}^{2}- \gamma^k_{11}\gamma_{k2}^{2} ) b^{1}_2$\\\\
$(22)\hspace{2em}( \gamma^2_{21}\gamma_{21}^{1}- \gamma^1_{22}\gamma_{11}^{2} ) b^{1}_1 +
		     ( \gamma^1_{12}\gamma_{12}^{2}- \gamma^2_{11}\gamma_{22}^{1} ) b^{2}_2 + $\\\\
$\hspace*{17em}    +( \gamma^k_{22}\gamma_{k1}^{1}- \gamma^k_{12}\gamma_{k2}^{1} ) b^{2}_1 +
		     ( \gamma^2_{2k}\gamma_{11}^{k}- \gamma^2_{1k}\gamma_{21}^{k} ) b^{1}_2\hfill$\\\\
$(12)\hspace{2em}2( \gamma^2_{2k}\gamma_{11}^{k}- \gamma^k_{21}\gamma_{k1}^{2} ) b^{1}_1 +
\hspace{5em}2( \gamma^2_{21}\gamma_{21}^{1}- \gamma^2_{11}\gamma_{22}^{1})b^{2}_1\hfill$\\\\
$(21)\hspace{8em}2( \gamma^1_{2k}\gamma_{21}^{k}- \gamma^1_{1k}\gamma_{22}^{k} ) b^{2}_2 + 
\hfill 2( \gamma^1_{22}\gamma_{11}^{2}- \gamma^1_{12}\gamma_{12}^{2})b^{1}_2$\\\\
Considered by itself, this system is degenerate.
Indeed, setting $$( \gamma^1_{22}\gamma_{11}^{2}- \gamma^1_{12}\gamma_{12}^{2} )=:A,
( \gamma^1_{2k}\gamma_{12}^{k}- \gamma^1_{1k}\gamma_{22}^{k} )=:B,
( \gamma^k_{21}\gamma_{k1}^{2}- \gamma^k_{11}\gamma_{k2}^{2} )=:C,$$
the system's matrix becomes:
$$
\left(
\begin{array}{cccc}

        A       &	        A            &          B         &          C         \\
       -A       &        -A		   &		  -B		&         -C         \\
       -2C      &		  0            &		  -2A		&          0         \\
        0       &  	  2B		   &		   0        &          2A        \\
\end{array} \right)
$$
The determinant of this is identically zero.

Let us now consider the other half of (\ref{eq:9}), the part $S(\Gamma_1) V_1$
( we will unite the halves afterward ).
Here we just adapt equation (\ref{eq:sys}) to the case $n=2,i=1,p=2$:\\\\
$(11)\hspace{2em}
(c_{11}^{12}-c_{21}^{11})b_{1}^{1}+
(c_{11}^{12}-c_{21}^{11})b_{2}^{2}+
(c_{12}^{12}-c_{22}^{11})b_{1}^{2}+
(c_{21}^{21}-c_{11}^{22})b_{2}^{1}
$\\\\
$(22)\hspace{2em}
(c_{12}^{22}-c_{22}^{21})b_{1}^{1}+
(c_{12}^{22}-c_{22}^{21})b_{2}^{2}+
(c_{22}^{11}-c_{12}^{12})b_{1}^{2}+
(c_{11}^{22}-c_{21}^{21})b_{2}^{1}
$\\\\
$(12)\hspace{1.5em}
2(c_{11}^{22}-c_{21}^{21})b_{1}^{1}+
\hspace{3.5em}[(c_{12}^{11}-c_{11}^{12})+
(c_{21}^{22}-c_{22}^{21})]b_{1}^{2}
\hspace{3.5em}$\\\\
$(21)\hspace{5em}
2(c_{12}^{12}-c_{22}^{11})b_{2}^{2}+
\hspace{3.5em} [(c_{22}^{21}-c_{21}^{22})+
(c_{11}^{12}-c_{12}^{11})]b_{2}^{1}
$\\\\
Setting 
$$c_{11}^{12}-c_{21}^{11}=:a, c_{12}^{12}-c_{22}^{11}=:b, c_{21}^{21}-c_{11}^{22}=:c,
c_{12}^{22}-c_{22}^{21}=:d\ ,$$ the above system's matrix can be written as:
\begin{equation}
\label{sys-matrix:gamma1}
\left(
\begin{array}{cccc}
        a       &	        a            &          b         &          c         \\
        d       &         d		   &		  -b		&         -c         \\
       -2c      &		  0            &		  d-a		&          0         \\
        0       &  	  2b		   &		   0        &         a-d        \\
\end{array} \right)
\end{equation}
It is also degenerate. The full (united) system of equations has the matrix:
$$\left(
\begin{array}{cccc}
        a+A       &	  a+A          &        b+B         &        c+C         \\
        d-A       &       d-A		   &	    -(b+B)		&      -(c+C)        \\
     -2(c+C)      &	  0            &     d-A-(a+A)	&          0         \\
        0         &  	2(b+B)	   &		0           &     a+A-(d-A)      \\
\end{array} \right)$$\\
It is also degenerate: in fact, it has exact same structure as (\ref{sys-matrix:gamma1}).
Even though not of full rank, generically the above system has rank 3, resulting in a 1-dimensional stabilizer.

Let us now consider the next, second jet of our connection. To calculate its stabilizer,
we need to solve the following equation from (\ref{sys:SYS}) for $V_4$ :
$${\mathcal{L}}_{V_{1}}\Gamma_{2}+\tilde{\mathcal{L}}_{V_{2}}\Gamma_{1}+
\tilde\mathcal{L}_{V_{3}}\Gamma_{0}+
{\displaystyle\frac{\partial^{2}V_{4}}{\partial x^{2}}}  =  0
$$
Its compatibility conditions are:
\begin{equation}
\label{eq:compV4}
({\mathcal{L}}_{V_{1}}\Gamma_{2}+\tilde{\mathcal{L}}_{V_{2}}\Gamma_{1}+
\tilde\mathcal{L}_{V_{3}}\Gamma_{0})_{ij,p}^l=  ({\mathcal{L}}_{V_{1}}\Gamma_{2}+\tilde{\mathcal{L}}_{V_{2}}\Gamma_{1}+
\tilde\mathcal{L}_{V_{3}}\Gamma_{0})_{pj,i}^l
\end{equation}
We will use the same strategy as in the previous section 
to prove that in this case stabilizer is trivial. Namely we will obtain a connection
2-jet, for which the above equation will be a non-degenerate homogeneous linear system.
We set $\Gamma_0=\Gamma_1=0$ , then (\ref{eq:compV4}) becomes:
\begin{equation}
\label{eq:compV4:restricted}
({\mathcal{L}}_{V_{1}}\Gamma_{2})_{ij,p}^l
=({\mathcal{L}}_{V_{1}}\Gamma_2)_{pj,i}^l
\end{equation}
We introduce notation for coefficients of $\Gamma_2$ : 
$${\Gamma_2}_{ij}^l=\sum_{s,t=1}^{2}d^l_{ijst}x^s x^t\ ,\ d^l_{ijst}=d^l_{jits}$$
With these, $$({\mathcal{L}}_{V_{1}}\Gamma_{2})_{ij,p}^l=
2d^l_{ijkt} b^k_p x^t+2d^l_{ijkp} b^k_t x^t-2d^k_{ijpt} b^l_k x^t+2d^l_{kjpt} b^k_i x^t+2d^l_{ikpt} b^k_j x^t$$ 
( $b_k^l$ are still coefficients of $V_1$, as in (\ref{V:comp}) ), and (\ref{eq:compV4:restricted}) 
( with $i=1,p=2$) is:\\\\
$(d^l_{kj2t}-d^l_{2jkt})b^k_1 + (d^l_{1jkt}-d^l_{kj1t})b^k_2 + (d^l_{1jk2}-d^l_{2jk1})b^k_t +$
\\\\
$\hspace*{17em}
(d^k_{2j1t}-d^k_{1j2t})b^l_k + (d^l_{1k2t}-d^l_{2k1t})b^k_j = 0 $\\\\
With the triple of indexes $(j,l,t)$ arbitrary, we have system of 8 equations in 4 variables: the coefficients of $V_1$. 
This is the system, equations are labelled by this index triple:\\\\
$
(111) 	
\hspace{2em}2(d^1_{1112}-d^1_{1211})b_1^1+(d^1_{1112}-d^1_{1211})b_2^2$\\
$\hspace*{17em}+(d^1_{1112}-d^1_{1211})b_1^2+  (d^2_{1211}-d^2_{1112})b_2^1=0$\\\\
$
(221)
\hspace{2em}2(d^2_{1212}-d^2_{2211})b_1^1+(d^2_{1212}-d^2_{2211})b_2^2+$\\
$\hspace*{12em}(d^2_{1222}-d^2_{2221}+d^1_{2211}-d^1_{1212})b_1^2+(d^2_{1112}-d^2_{1211})b_2^1=0
$\\\\
$
(121)
\hspace{2em}3(d^2_{1112}-d^2_{1211})b_1^1+\hspace{3em}(d^2_{1122}-d^2_{2211}+d^1_{2111}-d^1_{1112})b_1^2
\hspace{3em}=0
$\\\\
$
(211)
\hspace{2em}(d^1_{1212}-d^1_{2211})b_1^1+2(d^1_{1212}-d^1_{2211})b_2^2+$\\
$\hspace*{12em}(d^1_{1222}-d^1_{2221})b_1^2+ (d^1_{1112}-d^1_{2111}+d^2_{2211}-d^2_{1221})b_2^1=0
$\\\\
$
(112)
\hspace{2em}(d^1_{1122}-d^1_{1212})b_1^1+2(d^1_{1122}-d^1_{1212})b_2^2+$\\
$\hspace*{12em}(d^1_{1222}-d^1_{2221})b_1^2+(d^1_{1112}-d^1_{2111}+d^2_{1212}-d^2_{1122})b_2^1=0\\
$\\\\
$
(222)
\hspace{2em}(d^2_{1222}-d^2_{2221})b_1^1+2(d^2_{1222}-d^2_{2221})b_2^2+$\\
$\hspace*{17em}(d^1_{2221}-d^1_{1222})b_1^2+(d^2_{1122}-d^2_{2211})b_2^1=0\\
$\\\\
$
(122)
\hspace{2em}2(d^2_{1122}-d^2_{1212})b_1^1+(d^2_{1122}-d^2_{1212})b_2^2+$\\
$\hspace*{12em}(d^2_{1222}-d^2_{2221}+d^1_{1212}-d^1_{1122})b_1^2+(d^2_{1112}-d^2_{2111})b_2^1=0\\
$\\\\
$
(212)
\hspace{5.5em}3(d^1_{1222}-d^1_{2212})b_2^2+\hspace{4em}(d^2_{2221}-d^2_{1222}+d^1_{1122}-d^1_{2211})b_2^1=0$
\\\\
Setting:
$$
a=(d^1_{1112}-d^1_{1211}), e=(d^1_{1212}-d^1_{1122}), g=(d^1_{1222}-d^1_{2221}), h=(d^1_{1122}-d^1_{1212})\ ,$$
$$b=(d^2_{1211}-d^2_{1112}), c=(d^2_{1212}-d^2_{2211}), d=(d^2_{1222}-d^2_{2221}), f=(d^2_{1122}-d^2_{1212})
\ ,$$
we see the system take form:
$$\left(
\begin{array}{cccc}
2a & a & a   & b \\
2c & c & d+e &-b \\
-3b&   &f+c-a&   \\
-e & -2e & g & a-c\\
h  & 2h &  g & a-f\\
d  & 2d & -g & f+c\\
2f &  f & d-h& -b \\
   & 3g &    &h-e-d\\
\end{array} \right)$$
This is non-degenerate for a generic connection.
For example, if $a=d=c=e=f=h=0$, then it is:
$$\left(
\begin{array}{cccc}
 & & & b \\
& &  &-b \\
-3b&   &&   \\
 &  & g & \\
 &  &  g & \\
  &  & -g & \\
 &  & & -b \\
   & 3g &    &\\
\end{array} \right)$$
 
Now we can summarize what we know about exceptional stabilizers: 
\begin{prop}\label{except:stab:dim}
The stabilizer of a k-jet of a generic connections:\\
for $n=1$ is 1-dimensional ( equal to $G_1/G_2$ ) for any $k$ ;\\
for $n=2$ is equal to $G_1/G_2$ for $k=0$, is 1-dimensional for $k=1$ ,\\
and is trivial for $k\geq 2$ . 
\end{prop}
\section{Back to Poincar\'e series: proof of Theorem \ref{thm:series} }
\label{sec:Proof}
We will use Propositions \ref{stab:highdim} and \ref{except:stab:dim} to find dimension of a generic orbit:\\\\
$
\dim \mathcal{O}_{k}( \Gamma ) = \dim( T_{\mathrm{id}}(K_{k}/G_{\Gamma}))$\\\\
$\hspace*{5em}=\dim(\{ V | V={V}_{1}+{V}_{2}+\ldots+
{V}_{k+1}+{V}_{k+2} \}) - n^{2}\delta_{0}^{k}-\delta_{1}^{n}(1-\delta_{0}^{k})-\delta_{2}^{n}\delta_{1}^{k}\;,
$\\\\
where $\delta$ is a Kronecker symbol, taking care of non-zero stabilizers for various $k$ and $n$ ; 
${V_{i}}$ is an n-component vector, 
each component a homogeneous polynomial of degree $i$ in $x^{1},...,x^{n}$ .
So we have: 
$$
\dim\mathcal{O}_{k}=
{\displaystyle n\sum_{m=1}^{k+2} {n+m-1 \choose n-1} }-\delta_{1}^{n}-\delta_{2}^{n}\delta_{1}^{k} \qquad \textrm{for }k\geq1 \;, $$
$$  
\dim \mathcal{O}_{0}= 
n\sum_{m=1}^{2} {n+m-1 \choose n-1} - n^{2}=n{n+2-1 \choose n-1} =\frac{n^{2}(n+1)}{2}\ .$$
Dimension of the moduli space of connection $k$-jets $\mathcal{M}_{k}$ is:
$$\dim \mathcal{M}_{k}=\dim \mathcal{F}_{k}-\dim \mathcal{O}_{k}\;,$$  
where $\mathcal{F}_{k}$ is the space of connection k-jets.\\
For $k=0$:
$$\dim \mathcal{M}_{0}=\dim \mathcal{F}_{0}-\dim \mathcal{O}_{0}=n\frac{n(n+1)}{2}-\frac{n^{2}(n+1)}{2}=0\ .$$
For $k\geq1$ :
$$\dim \mathcal{M}_{k}=\dim \mathcal{F}_{k}-\dim \mathcal{O}_{k}$$
$$=n\frac{n(n+1)}{2}\sum_{m=0}^{k}{n+m-1 \choose n-1}-n\sum_{m=1}^{k+2}{n+m-1 \choose n-1}+
\delta_{1}^{n}+\delta_{2}^{n}\delta_{1}^{k}\ .$$
The Poincar\'e series is:
$$p_{\Gamma}(t)=\dim \mathcal{M}_{0}+\sum_{k=1}^{\infty}(\dim \mathcal{M}_{k}-\dim \mathcal{M}_{k-1})t^{k}$$
\begin{equation}
\label{eq:P'series}
=n\sum_{k=1}^{\infty}\;\biggl[\frac{n(n+1)}{2}{n+k-1 \choose n-1}
-n{n+k+1 \choose n-1}\biggr]t^{k}+(t-t^2)\delta_{2}^{n}\;,
\end{equation}
as required.
Simplifying (\ref{eq:P'series}), we obtain the following\\\\
{\bf Fact}
\begin{em}
The Poncar\'e series $p_\Gamma(t)$ is a rational function.
Namely,
\begin{equation}
\label{P'ratio}
p_{\Gamma}(t)=(t-t^2)\delta_{2}^{n}+nD_{\Gamma}\left(\frac{1}{1-t}\right)+\frac{(n-1)n^{2}(n+1)}{2}
\end{equation}
where $D_{\Gamma}$ is a differential operator of order $n-1$ :
$$D_{\Gamma}=\frac{n(n+1)}{2}{n+t\frac{d}{dt}-1 \choose n-1}
-n{n+t\frac{d}{dt}+1 \choose n-1}\ ,$$
with $${n+t\frac{d}{dt}-1 \choose n-1}=\frac{1}{(n-1)!}(t\frac{d}{dt}+1)\ldots(t\frac{d}{dt}+n-1)\ ,$$ 
$${n+t\frac{d}{dt}+1 \choose n-1}=\frac{1}{(n-1)!}(t\frac{d}{dt}+3)\ldots(t\frac{d}{dt}+n+1) \ .$$
\end{em}
Indeed, denote 
$$\varphi_{m}(t)=\sum_{k=0}^{\infty}k^{m}t^{k}\ ,\qquad m\in\mathbb{Z}_{+}\ ,$$
then
$$\varphi_{m}(t)=\sum_{k=0}^{\infty}k^{m-1}kt^{k-1}t=t\left(\sum_{k=0}^{\infty}k^{m-1}t^{k}\right)'
=\left(t\frac{d}{dt}\right)\varphi_{m-1}(t)\quad \mathrm{for}\ m\in\mathbb{N}\,.$$
Thus
$$\varphi_{m}(t)=\left(t\frac{d}{dt}\right)^{m}\varphi_{0}(t)=
\left(t\frac{d}{dt}\right)^{m}\left(\frac{1}{1-t}\right)\,.$$
Hence,
$$\sum_{k=0}^{\infty}\;\biggl[\frac{n(n+1)}{2}{n+k-1 \choose n-1}-n{n+k+1 \choose n-1}\biggr]t^{k}$$
$$=\biggl[\frac{n(n+1)}{2}{n+t\frac{d}{dt}-1 \choose n-1}-n{n+t\frac{d}{dt}+1 \choose n-1}\biggr]
\left(\frac{1}{1-t}\right)\,.$$
We have: $$a_{0}=\frac{n(n+1)}{2}-n{n+1 \choose n-1}=-\frac{(n-1)n(n+1)}{2}\,.$$
So $$p_{\Gamma}(t)=(t-t^2)\delta_{2}^{n}+
n\sum_{k=0}^{\infty}\;\biggl[\frac{n(n+1)}{2}{n+k-1 \choose n-1}-n{n+k+1 \choose n-1}\biggr]t^{k}
+\frac{(n-1)n^{2}(n+1)}{2}$$
$$=(t-t^2)\delta_{2}^{n}+D_{\Gamma}\left(\frac{1}{1-t}\right)+\frac{(n-1)n^{2}(n+1)}{2}$$
\hfill$\Box$

\end{document}